\documentclass{amsart}

\usepackage{fancyhdr}
\usepackage{lastpage}
\usepackage{stmaryrd,yhmath}

\newcommand{\bdis}{\begin{displaymath}}
\newcommand{\edis}{\end{displaymath}}
\newcommand{\be}{\begin{equation}}
\newcommand{\ee}{\end{equation}}
\newcommand{\mbb}{\mathbb}
\newcommand{\mcal}{\mathcal}

\newcommand{\vp}{\varphi}

\newcommand{\vth}{\vartheta}

\theoremstyle{definition}

\theoremstyle{remark}
\newtheorem{remark}[]{Remark}

\newtheorem*{mydef11}{{\bf Theorem 1}}

\newtheorem*{mydef12}{{\bf Theorem 2}}

\newtheorem*{mydef13}{{\bf Theorem 3}}

\newtheorem*{mydef14}{{\bf Theorem 4}}

\newtheorem*{mydef41}{{\bf Corollary 1}}

\newtheorem*{mydef42}{{\bf Corollary 2}}

\numberwithin{equation}{section}

\begin{document}

\title{Properties of the sequence $\{ Z[t_\nu(\tau)]\}$, Jacob's ladders and new kind of infinite set of metamorphosis of
main multiform}

\author{Jan Moser}

\address{Department of Mathematical Analysis and Numerical Mathematics, Comenius University, Mlynska Dolina M105, 842 48 Bratislava, SLOVAKIA}

\email{jan.mozer@fmph.uniba.sk}


\begin{abstract}
In this paper we study properties of some sums of members of the sequence $\{ Z[t_\nu(\tau)]\}$. Our results are expressed in statements
proving essential influence of the Lindel\" of hypothesis on corresponding formulae. In this paper: the parts 1 -- 6 are English
version of our paper \cite{6}, and the part 7 of this work contains current results, namely new set of metamorphosis of the main
multiform from our paper \cite{7}.
\end{abstract}
\maketitle

\section{Introduction}

\subsection{}

Let us remind the Riemann-Siegel formula
\be \label{1.1} \begin{split}
& Z(t)=2\sum_{n\leq\rho(t)}\frac{1}{\sqrt{n}}\cos\{\vth(t)-t\ln n\}+\mcal{O}(t^{-1/4}), \\
& \vth(t)=-\frac t2\ln\pi+\mbox{Im}\-\ln\Gamma\left(\frac 14+i\frac t2\right),\ \rho(t)=\sqrt{\frac{t}{2\pi}},
\end{split}
\ee
(see \cite{8}, pp. 79, 329). In the paper \cite{5} we have defined the following two notions in connection
with the formula (\ref{1.1}):

\begin{itemize}
\item[(a)] the class of sequences
\bdis
\{ t_\nu(\tau)\}
\edis
by the condition
\be \label{1.2}
\begin{split}
& \vth[t_\nu(\tau)]=\pi\nu+\tau,\ \tau\in [-\pi,\pi],\ t(0)=t_\nu, \\
& \nu_0\leq \nu,\ \nu_0\in \mbb{N},
\end{split}
\ee
where $\nu_0$ is sufficiently big and fixed number,

\item[(b)] the following classes of disconnected sets
\be \label{1.3}
\begin{split}
& \mbb{G}_1=\mbb{G}_1(x,T,H)=\\
& =\bigcup_{T\leq t_{2\nu}\leq T+H} \{ t:\ t_{2\nu}(-x)<t<t_{2\nu}(x)\},\ x\in (0,\pi/2], \\
& \mbb{G}_2=\mbb{G}_1(y,T,H)=\\
& =\bigcup_{T\leq t_{2\nu+1}\leq T+H} \{ t:\ t_{2\nu+1}(-y)<t<t_{2\nu+1}(y)\},\ y\in (0,\pi/2], \\
& H\in (0,H_1],\ H_1=T^{1/6+\epsilon},
\end{split}
\ee
where $\epsilon$ is positive and arbitrarily small number.
\end{itemize}

Further, by making use of a synthesis of discrete and continuous method, that is (see (\ref{1.2})):
\begin{itemize}
\item[(c)] summation according to $\nu$,
\item[(d)] integration according to $\tau$,
\end{itemize}
we have obtained the following new mean-value formulae (see \cite{5})
\be \label{1.4}
\begin{split}
& \frac{1}{m\{\mbb{G}_1\}}\int_{\mbb{G}_1}Z(t){\rm d}t\sim 2\frac{\sin x}{x}, \\
& \frac{1}{m\{\mbb{G}_2\}}\int_{\mbb{G}_2}Z(t){\rm d}t\sim -2\frac{\sin y}{y}, \\
& T\to\infty,
\end{split}
\ee
on discrete sets $\mbb{G}_1,\mbb{G}_2$, where $m\{ \cdot \}$ stands for the measure of the corresponding set.

\begin{remark}
Formulae (\ref{1.4}) are the first mean-value theorems linear in the function $Z(t)$ in the theory of the Riemann zeta-function.
\end{remark}

\subsection{}

In this paper we shall study the following sums
\be \label{1.5}
\begin{split}
& \sum_{T\leq t_\nu\leq T+H}\{ Z[t_\nu(\tau)]-Z(t_\nu)\}, \\
& \sum_{T\leq t_\nu\leq T+H} (-1)^\nu \{ Z[t_\nu(\tau)]-Z(t_\nu)\}, \\
& H\in (0,T^{1/6+\epsilon}],
\end{split}
\ee
and corresponding integral formulae also.

\begin{remark}
The Lindel\" of hypothesis gives an essential influence on resulting asymptotic formulae.
\end{remark}

\subsection{}

It is the most astonishing fact for the author that the expression
\be \label{1.6}
Z[t_\nu(\tau)]-Z(t_\nu)
\ee
coming from our old (1983) paper \cite{6} is close connected with the $\zeta$-factorization and with the
metamorphosis of main multiform from our last (2015) paper \cite{7}.

\begin{remark}
Namely, from (\ref{1.6}) it follows that there is an essentially new set of metamorphosis of the main multiform (see \cite{7}, (2.4)).
\end{remark}

\section{Theorem 1}

Let
\be \label{2.1}
F(\tau,T,H)=\sum_{T\leq t_\nu\leq T+H} Z[t_\nu(\tau)]
\ee
and
\be \label{2.2}
S(a,b)=\sum_{1\leq a\leq n<b\leq 2a}n^{it},\ b\leq\sqrt{\frac{t}{2\pi}}
\ee
stands for elementary trigonometric sum. The following theorem holds true.

\begin{mydef11}
If
\be \label{2.3}
|S(a,b)|<A(\Delta)\sqrt{a}t^\Delta,\ \Delta\in (0,1/6] ,
\ee
(see (\ref{2.2})), then
\be \label{2.4}
\begin{split}
 & F(\tau,T,H)-F(0,T,H)=\mcal{O}(T^\Delta\ln T), \\
 & H\in (0,T^{1/6+\epsilon}]
\end{split}
\ee
uniformly for $\tau\in [-\pi,\pi]$.
\end{mydef11}

In the case $\Delta=1/6$ we have the following

\begin{mydef41}
\bdis
F(\tau,T,H)-F(0,T,H)=\mcal{O}(T^{1/6+\epsilon})
\edis
uniformly for $\tau\in [-\pi,\pi]$.
\end{mydef41}

If the Lindel\" of hypothesis holds true, then
\bdis
\Delta=\frac{\epsilon}{2},
\edis
(comp. \cite{1}, p. 89), and we have the following

\begin{mydef42}
On Lindel\" of hypothesis
\bdis
F(\tau,T,H)-F(0,T,H)=\mcal{O}(T^{\epsilon})
\edis
uniformly for $\tau\in [-\pi,\pi]$.
\end{mydef42}

\section{Theorem 2}

First of all the following formula holds true
\be \label{3.1}
\sum_{T\leq t_\nu\leq T+H}(-1)^\nu Z(t_\nu)=\frac 1\pi H\ln\frac{T}{2\pi}+\mcal{O}(T^{\Delta}\ln T)
\ee
in the case (\ref{2.3}) (see \cite{3}, p. 89, (5)). Next, we have proved (see \cite{5}, (5.1)) that in the
case (\ref{2.3}) the formula
\be \label{3.2}
\begin{split}
 & \sum_{T\leq t_\nu\leq T+H}(-1)^\nu Z[t_\nu(\tau)]= \\
 & = \frac 1\pi H\ln\frac{T}{2\pi}\cos\tau+\mcal{O}(T^{\Delta}\ln T)
\end{split}
\ee
follows, where the $\mcal{O}$-estimate is true uniformly for $\tau\in [-\pi,\pi]$.

Further, from (\ref{3.1}), (\ref{3.2}) the formula
\be \label{3.3}
\begin{split}
& \sum_{T\leq t_\nu\leq T+H}(-1)^\nu \{ Z[t_\nu(\tau)]-Z(t_\nu)\}= \\
& = -\frac{2}{\pi} H\ln\frac{T}{2\pi}\sin^2\frac{\tau}{2}+\mcal{O}(T^{\Delta}\ln T)
\end{split}
\ee
follows. Consequently, we have by (\ref{2.4}), (\ref{3.3}) the following

\begin{mydef12}
It follows from (\ref{2.3}) that
\be \label{3.4}
\begin{split}
& \sum_{T\leq t_{2\nu}\leq T+H}(-1)^\nu \{ Z[t_{2\nu}(\tau)]-Z(t_{2\nu})\}= \\
& = -\frac 1\pi H\ln\frac{T}{2\pi}\sin^2\frac{\tau}{2}+\mcal{O}(T^{\Delta}\ln T), \\
& \sum_{T\leq t_{2\nu+1}\leq T+H}(-1)^\nu \{ Z[t_{2\nu+1}(\tau)]-Z(t_{2\nu+1})\}= \\
& = \frac 1\pi H\ln\frac{T}{2\pi}\sin^2\frac{\tau}{2}+\mcal{O}(T^{\Delta}\ln T),
\end{split}
\ee
where the $\mcal{O}$-estimates hold true uniformly for $\tau\in [-\pi,\pi]$.
\end{mydef12}

\begin{remark}
Our formulae (\ref{3.3}), (\ref{3.4}) are asymptotic ones in the case
\be \label{3.5}
H=T^\Delta \ln T.
\ee
\end{remark}

\section{Theorem 3}

Of course,
\be \label{4.1}
\begin{split}
& \frac{1}{t_{2\nu}(x)-t_{2\nu}(-x)}\int_{t_{2\nu}(-x)}^{t_{2\nu}(x)} Z(t){\rm d}t=Z[\xi_{2\nu}(x)], \\
& \frac{1}{t_{2\nu+1}(y)-t_{2\nu+1}(-y)}\int_{t_{2\nu+1}(-y)}^{t_{2\nu+1}(y)} Z(t){\rm d}t=Z[\xi_{2\nu+1}(y)],
\end{split}
\ee
where
\bdis
\begin{split}
& \xi_{2\nu}(x)\in (t_{2\nu}(-x),t_{2\nu}(x)), \\
& \xi_{2\nu+1}(y)\in (t_{2\nu+1}(-y),t_{2\nu+1}(y)),
\end{split}
\edis
and the numbers
\bdis
Z[\xi_{2\nu}(x)],\ Z[\xi_{2\nu+1}(y)]
\edis
are the mean-values of the function $Z(t)$ with respect to corresponding segments. Also, in this direction, there are asymptotic formulae
for the sums of differences
\be \label{4.2}
\begin{split}
& Z[\xi_{2\nu}(x)]-Z(t_{2\nu}), \\
& Z[\xi_{2\nu+1}(y)]-Z(t_{2\nu+1}).
\end{split}
\ee

\begin{remark}
The behavior of the differences (\ref{4.2}) is very irregular.
\end{remark}

However, the following theorem holds true.

\begin{mydef13}
It follows from (\ref{2.3}) that
\be \label{4.3}
\begin{split}
& \sum_{T\leq t_{2\nu}\leq T+H}\{ Z[\xi_{2\nu}(x)]-Z(t_{2\nu})\}= \\
& = -\frac{1}{2\pi}\left( 1-\frac{\sin x}{x}\right) H\ln\frac{T}{2\pi}+\mcal{O}(T^\Delta\ln T), \\
& \sum_{T\leq t_{2\nu+1}\leq T+H}\{ Z[\xi_{2\nu+1}(y)]-Z(t_{2\nu+1})\}= \\
& = \frac{1}{2\pi}\left( 1-\frac{\sin y}{y}\right) H\ln\frac{T}{2\pi}+\mcal{O}(T^\Delta\ln T).
\end{split}
\ee
\end{mydef13}

\section{Proof of Theorem 1}

First of all, we obtain from (\ref{1.1}) (see \cite{2}, (57)) that
\bdis
\begin{split}
& Z(t)=2\sum_{n<P_0} \frac{1}{\sqrt{n}}\cos\{\vth(t)-t\ln n\}+\mcal{O}(t^{-1/4}), \\
& t\in [T,T+H],\ H\in (0,\sqrt[4]{T}],\ P_0=\sqrt{\frac{T}{2\pi}}.
\end{split}
\edis
Since (see (\ref{1.2}))
\bdis
\vth[t_\nu(\tau)]-\vth(t_\nu)=\tau,
\edis
then (see \cite{2}, (41), (42))
\be \label{5.1}
t_\nu(\tau)-t_\nu=\frac{\tau}{\ln P_0}+\mcal{O}\left(\frac{H}{T\ln^2T}\right),\ t_\nu\in [T,T+H].
\ee
Next,
\bdis
\begin{split}
& \sin\left(\frac{\pi}{2}-\frac{t_\nu(\tau)-t_\nu}{2}\ln n\right)=\sin\left\{ \frac{\tau}{\pi}X(n)\right\}+\mcal{O}\left(\frac{H}{T\ln T}\right), \\
& \sin\left(\frac{\pi}{2}-\frac{t_\nu(\tau)+t_\nu}{2}\ln n\right)=\sin\left\{ \frac{\tau}{\pi}X(n)-t_\nu\ln n\right\}+\mcal{O}\left(\frac{H}{T\ln T}\right), \\
& X(n)=\frac{\pi}{2\ln P_0}\ln\frac{P_0}{n}, \\
& 0<X(n)\leq \frac{\pi}{2},\ 1\leq n<P_0.
\end{split}
\edis
Now (see (\ref{1.2}))
\be \label{5.2}
\begin{split}
& Z[t_\nu(\tau)]-Z(t_\nu)= \\
& = 4(-1)^{\nu+1}\sum_{n<P_0}\frac{1}{\sqrt{n}}\sin\left\{ \frac{\tau}{\pi}X(n)\right\}\sin\left\{ \frac{\tau}{\pi}X(n)-t_\nu\ln n\right\}+ \\
& + \mcal{O}\left(\sqrt[4]{T}\frac{H}{T\ln T}\right)+\mcal{O}(T^{-1/4})= \\
& =4(-1)^{\nu+1}\sum_{n<P_0}\frac{1}{\sqrt{n}}\sin^2\left\{ \frac{\tau}{\pi}X(n)\right\}\cos\{ t_\nu\ln n\}+ \\
& + 2(-1)^\nu\sum_{n<P_0}\frac{1}{\sqrt{n}}\sin\left\{ \frac{2\tau}{\pi}X(n)\right\}\sin\{ t_\nu\ln n\}+ \\
& + \mcal{O}(T^{-1/4}).
\end{split}
\ee
Hence (see (\ref{2.1}), \cite{2}, (59))
\bdis
\begin{split}
& F(\tau,T,H)-F(0,T,H)=\sum_{T\leq t_\nu\leq T+H}\{ Z[t_\nu(\tau)]-Z(t_\nu)\}= \\
& = -4\sum_{n<P_0}\frac{1}{\sqrt{n}}\sin^2\left\{ \frac{\tau}{\pi}X(n)\right\}\cdot\sum_{T\leq t_\nu\leq T+H}(-1)^\nu \cos\{ t_\nu\ln n\}+ \\
& + 2\sum_{n<P_0} \sin\left\{ \frac{\tau}{\pi}X(n)\right\}\cdot \sum_{T\leq t_\nu\leq T+H}(-1)^\nu \sin\{ t_\nu\ln n\}+ \\
& + \mcal{O}(\ln T)= \\
& = -4w_1+2w_2+\mcal{O}(\ln T).
\end{split}
\edis
The sum $w_1$ contains the following typical member
\bdis
w_{11}=\sum_{n<P_0}\sin^2\left\{ \frac{\tau}{\pi}X(n)\right\}\frac{\tan\frac{\omega}{2}}{\sqrt{n}}\sin\vp,
\edis
(comp. \cite{2}, (54)), where
\bdis
\begin{split}
& \vp=t_{\bar{\nu}}\ln n, \\
& \frac{\omega}{2}=\frac{\pi}{2}\frac{\ln n}{\ln P_0}=\frac{\pi}{2}-X(n), \\
& \tan\frac{\omega}{2}=\cot X(n),
\end{split}
\edis
(comp. \cite{2}, (43), (50)). Consequently,,
\be \label{5.3}
w_{11}=\left(\frac{\tau}{\pi}\right)^2\sum_{n<P_0} X\frac{\sin^2\left(\frac{\tau}{\pi}X\right)}{\left(\frac{\tau}{\pi}X\right)^2}
\frac{X}{\sin X}\cos X\frac{1}{\sqrt{n}}\sin\vp.
\ee
Since in the case (\ref{2.3}) we have that
\bdis
\sum_{1\leq n<P_1\leq P_0}\frac{1}{\sqrt{n}}\sin\vp=\mcal{O}(T^\Delta\ln T)
\edis
then, by making use the Abel's transformation several-time in (\ref{5.3}), we obtain the estimate
\bdis
w_{11}=\mcal{O}(T^\Delta\ln T),
\edis
and, consequently,
\bdis
w_{1}=\mcal{O}(T^\Delta\ln T)
\edis
uniformly for $\tau\in [-\pi,\pi]$. Because of the simple identity
\bdis
\sin(t_\nu\ln n)=\cos(t_\nu\ln n-\pi/2)
\edis
we obtain by a similar way the estimate for $w_2$ too.

\section{Proof of Theorem 3}

Since (see \cite{5}, (2.1))
\bdis
t_{2\nu}(x)-t_{2\nu}(-x)=\frac{2x}{\ln P_0}+\mcal{O}\left(\frac{xH}{T\ln^2T}\right),
\edis
and by (\ref{2.3})
\bdis
Z(t)=\mcal{O}(T^\Delta\ln T) ,
\edis
then (see (\ref{2.3}), (\ref{4.1}), \cite{5}, (6.2))
\bdis \begin{split}
& \int_{-x}^x Z[t_{2\nu}(\tau)]{\rm d}\tau=\ln P_0\int_{t_{2\nu}(-x)}^{t_{2\nu}(x)}Z(t){\rm d}t+
\mcal{O}(xHT^{-5/6})= \\
& = [t_{2\nu}(x)-t_{2\nu}(-x)]Z[\xi_{2\nu}(x)]\ln P_0+\mcal{O}(xHT^{-5/6})= \\
& = 2xZ[\xi_{2\nu}(x)]+\mcal{O}(xHT^{-5/6}).
\end{split}
\edis
Hence,
\bdis
\begin{split}
 & \int_{-x}^x \sum_{T\leq t_{2\nu}\leq T+H}\{ Z[t_{2\nu}(\tau)]-Z(t_{2\nu})\}{\rm d}\tau=\\
 & = 2x\sum_{T\leq t_{2\nu}\leq T+H}\{ Z[\xi_{2\nu}(x)]-Z(t_{2\nu})\}+ \\
 & + \mcal{O}(xHT^{-5/6}H\ln T),
\end{split}
\edis
where we have used the formula (see \cite{2}, (2.3))
\bdis
\sum_{T\leq t_{\nu}\leq T+H}1\sim \frac{1}{2\pi} H\ln\frac{T}{2\pi}.
\edis
Consequently, the integration of the first formula in (\ref{3.4}) by $\tau\in [-\pi,\pi]$ together
with the formula
\bdis
\int_{-x}^x \sin^2\frac{\tau}{2}{\rm d}\tau=x-\sin x
\edis
gives the first formula in (\ref{4.3}), and, by a similar way, we obtain also the second formula in (\ref{4.3}).

\section{On new kind of metamorphosis of main multiform}

\subsection{}

First of all, we define the following sets
\be \label{7.1}
\begin{split}
 & w_\nu(T,H)=\{\tau:\ \tau\in (0,\pi],\ Z[t_\nu(\tau)]\not= Z(t_\nu)\}, \\
 & t_\nu\in [T,T+H],\ H\in\left[\frac{\ln\ln T}{\ln T},T^{\frac{1}{\ln\ln T}}\right].
\end{split}
\ee
The lower bound for $H$ in (\ref{7.1}) is chosen with respect to the formula (see (\ref{5.1})
\be \label{7.2}
t_\nu(\tau)-t_\nu\sim\frac{\tau}{\ln P_0},\ T\to\infty .
\ee
Next, by the Newton-Leibniz formula
\be \label{7.3}
\left|\int_{t_\nu}^{t_\nu(\tau)}Z'(u){\rm d}u\right|=|Z[t_\nu(\tau)]-Z(t_\nu)|,\ \tau\in w_\nu.
\ee

\begin{remark}
It is a new thing that the elementary identity (\ref{7.3}) has a nontrivial continuation by means of
Jacob's ladders and reversely iterated integrals.
\end{remark}

\subsection{}

Since (see (\ref{7.2}))
\bdis
t_\nu(\tau)-t_\nu=o\left(\frac{T}{\ln T}\right),
\edis
then we obtain from (\ref{7.3}) by \cite{7}, (5.2), (5.3) that
\be \label{7.4}
\begin{split}
 & \prod_{r=0}^{k-1} |Z[\vp_1^r(\beta)]|\sim \\
 & \sim \sqrt{\frac{|Z[t_\nu(\tau)]-Z[t_\nu]|}{\overset{k}{\wideparen{t_\nu(\tau)}}-\overset{k}{\wideparen{t_\nu}}}}
 \ln^kT\frac{1}{\sqrt{|Z'[\vp_1^k(\beta)]|}},\ \tau\in w_\nu.
\end{split}
\ee
Next, we obtain from (\ref{7.4}) by \cite{7}, (5.4) -- (5.7) that
\bdis
\begin{split}
 & \prod_{l=1}^k |Z[\bar{\alpha}_l]|\sim \frac{\Lambda_1}{\sqrt{|Z'(\bar{\alpha}_0)|}},\ \tau\in w_\nu, \\
 & \bar{\alpha}_r=\bar{\alpha}_r(T,H,k,\nu,\tau),\ r=0,1,\dots,k, \\
 & \bar{\alpha}_r\not=\gamma:\ \zeta\left(\frac 12+i\gamma\right)=0, \\
 & \Lambda_1=\Lambda_1(T,H,k,\nu,\tau)=
 \sqrt{\frac{|Z[t_\nu(\tau)]-Z[t_\nu]|}{\overset{k}{\wideparen{t_\nu(\tau)}}-\overset{k}{\wideparen{t_\nu}}}}
 \ln^kT .
\end{split}
\edis
Consequently, the following theorem holds true (comp. \cite{7}, Theorem)

\begin{mydef14}
Let
\bdis
k=1,\dots,k_0,\ k_0\in\mbb{N}
\edis
for every fixed $k_0$ and let
\bdis
H=H(T)\in\left[\frac{\ln\ln T}{\ln T},T^{\frac{1}{\ln\ln T}}\right]
\edis
for every sufficiently big $T>0$. Then for every $k,H$ and for every $\tau\in w_\nu$ there are functions
\bdis
\begin{split}
 & \bar{\alpha}_r=\bar{\alpha}_r(T,H,k,\nu,\tau)>0,\ r=0,1,\dots,k , \\
 & \bar{\alpha}_r\not=\gamma:\ \zeta\left(\frac 12+i\gamma\right)=0,
\end{split}
\edis
such that the following $\zeta$-factorization formula
\bdis
\begin{split}
 & \frac{\Lambda_1}{\sqrt{|Z'(\bar{\alpha}_0)|}}\sim \prod_{l=1}^k |Z(\bar{\alpha}_l)|,\ \tau\in w_\nu, \\
 & \Lambda_1=\Lambda_1(T,H,k,\nu,\tau)=
 \sqrt{\frac{|Z[t_\nu(\tau)]-Z[t_\nu]|}{\overset{k}{\wideparen{t_\nu(\tau)}}-\overset{k}{\wideparen{t_\nu}}}}
 \ln^kT
\end{split}
\edis
holds true. Moreover, the sequence
\bdis
\{\bar{\alpha}_r\}_{r=0}^k
\edis
obeys the following properties
\bdis
\begin{split}
 & T<\bar{\alpha}_0<\bar{\alpha}_1<\dots<\bar{\alpha}_k, \\
 & \bar{\alpha}_{r+1}-\bar{\alpha}_r\sim (1-c)\pi(T),\ r=0,1,\dots,k-1
\end{split}
\edis
where
\bdis
\pi(T)\sim \frac{T}{\ln T},\ T\to\infty
\edis
is the prime-counting function and $c$ is the Euler's constant.
\end{mydef14}

\subsection{}

Let us remind that we have proved in our paper \cite{5} the following formula of the Riemann-Siegel type
\bdis
\begin{split}
 & Z'(t)=-2\sum_{n\leq\rho(t)}\frac{1}{\sqrt{n}}\ln\frac{\rho(t)}{n}\sin\{\vth(t)-t\ln n\}+ \\
 & + \mcal{O}(t^{-1/4}\ln t),\ \rho(t)=\sqrt{\frac{t}{2\pi}}.
\end{split}
\edis
Let (comp. \cite{7}, (2.1) -- (2.3))
\bdis
\begin{split}
 & a_n=\frac{2}{a_n} , \\
 & \bar{f}_n(t)=\ln\frac{\rho(t)}{n}\sin\{\vth(t)-t\ln n\}, \\
 & \bar{R}(t)=\mcal{O}(t^{-1/4}\ln t).
\end{split}
\edis
Then we obtain (comp. \cite{7}, (2.5)) from Theorem 4 the following

\begin{mydef42}
\be \label{7.5}
\begin{split}
 & \prod_{r=1}^k \left|\sum_{n\leq \rho(\bar{\alpha}_r)}a_nf_n(\bar{\alpha}_r)+R(\bar{\alpha}_r)\right|\sim \\
 & \sim \frac{\Lambda_1}
 {\sqrt{\left|\sum_{n\leq \rho(\bar{\alpha}_0)}a_n\bar{f}_n(\bar{\alpha}_0)+\bar{R}(\bar{\alpha}_0)\right|}},
\end{split}
\ee
i.e. an infinite set of metamorphoses of the main multiform
(comp. \cite{7}, (2.4))
\bdis
G(x_1,\dots,x_k)=\prod_{r=1}^k |Z(x_r)|,\ \bar{x}_r>T>0,\ k\geq 2
\edis
into quite distinct monoform on the right-hand side of (\ref{7.5}) corresponds to the infinite subset of the points
\bdis
\{ \bar{\alpha}_1(T),\bar{\alpha}_2(T),\dots,\bar{\alpha}_k(T)\},\ T\in (T_0,+\infty) ,
\edis
where $T_0$ is sufficiently big.
\end{mydef42}

\thanks{I would like to thank Michal Demetrian for his help with electronic version of this paper.}


\begin{thebibliography}{29}
\bibitem{1}
A.A, Karatsuba,
`\emph{Basic analytic number theory}`, Moscow, 1975 (in Russian).
%
\bibitem{2}
J. Moser,
`On one sum in the theory of the Riemann zeta-function`,  Acta Arith., 31 (1976), 31-43, (in Russian).
Acta Math. 41 (1918), 119-196.
%
\bibitem{3}
J. Moser,
`On one theorem of Hardy-Littlewood in the theory of the Riemann zeta-function`,
Acta Arith., 31 (1976), 45-51; 40 (1981), 97-107, (in Russian).
%
\bibitem{4}
J. Moser, `On roots of the equation $Z'(t)=0$`, Acta Arith., 40 (1981), 79-89, arXiv: 1303.0967.
%
\bibitem{5}
J. Moser, `New consequences of the Riemann-Siegel formula and a law of asymptotic equality of signum areas`,
Acta Arith., 42 (1982), 1-10, arXiv: 1312.4767.
%
\bibitem{6}
J. Moser, `On properties of the sequence $\{ Z[t_\nu(\tau)]\}$ in the theory of the Riemann zeta-function`,
Acta Math. Univ. Comen., 42-43 (1983), 55-63, (in Russian).
%
\bibitem{7}
J. Moser, `Jacob's ladders, $\zeta$-factorization and infinite set of metamorphosis of a multiform`,
arXiv: 1501.07705.
%
\bibitem{8}
E.C. Titchmarsh, `\emph{The theory of the Riemann zeta-function}`, Clarendon Press, Oxford, 1951.



\end{thebibliography}
\end{document}